\documentclass{gSTA2e}

\usepackage{epstopdf}
\usepackage{subfigure}

\theoremstyle{plain}
\newtheorem{theorem}{Theorem}[section]

\newtheorem{lemma}[theorem]{Lemma}

\theoremstyle{definition}

\theoremstyle{remark}

\begin{document}



\title{Martingale inequalities of type Dzhaparidze  and van Zanten}

\author{
\name{X. Fan\textsuperscript{a}$^{\ast}$\thanks{$^\ast$Corresponding author. Email: fanxiequan@hotmail.com},
I. Grama\textsuperscript{b}
and Q. Liu\textsuperscript{b}}
\affil{\textsuperscript{a}Regularity Team, Inria, France;
\textsuperscript{b}Universit\'{e} de Bretagne-Sud,  UMR 6205, LMBA, Vannes, France.}
\received{v5.0 released February 2015}
}

\maketitle

\begin{abstract}
Freedman's inequality is a supermartingale counterpart to Bennett's inequality. This result shows that
the  tail probabilities of a supermartingale is controlled by the quadratic characteristic
and a uniform upper bound for the supermartingale difference sequence. Replacing the quadratic characteristic  by $\textrm{H}_k^y:= \sum_{i=1}^k\left(\mathbf{E}(\xi_i^2 |\mathcal{F}_{i-1}) +\xi_i^2\textbf{1}_{\{|\xi_i|> y\}}\right),$
Dzhaparidze and van Zanten  (\emph{Stochastic Process.  Appl.}, 2001) have extended Freedman's inequality to martingales with unbounded differences.  In this paper, we prove that  $\textrm{H}_k^y$ can be refined to $\textrm{G}_k^{y}   :=\sum_{i=1}^k   \left( \mathbf{E}(\xi_i^2\textbf{1}_{\{\xi_i \leq y\}} |\mathcal{F}_{i-1}) + \xi_i^2\textbf{1}_{\{\xi_i> y\}}\right).$   Moreover, we also establish two inequalities of type Dzhaparidze and van Zanten. These results  extend   Sason's  inequality (\emph{Statist. Probab. Lett.}, 2012)  to the martingales with possibly unbounded differences and establish  the connection between  Sason's  inequality and  De la Pe\~{n}a's inequality (\emph{Ann.\ Probab.,} 1999).  An application to self-normalized deviations is given.
\end{abstract}

\begin{keywords}
Freedman's inequality; De la Pe\~{n}a's inequality; exponential inequalities; tail probabilities; martingales; self-normalized deviations
\end{keywords}

\begin{classcode}  60G42; 60E15 
\end{classcode}

\section{Introduction}
Exponential inequalities for tail probabilities of sums of independent real-valued  random variables
and their extension to martingales have numerous important applications in probability
theory and statistics. See, for instance,  De la Pe\~{n}a and Pang \cite{D09},  Bercu  and Touati \cite{BT08} and \cite{F14}.
The classical   Bennett inequality  \cite{B62} gives a tail bound for sums of independent random
variables with a bounded range.
If $ (\xi _i)_{i \geq1}$ are zero-mean independent random
variables, all bounded  by some constant $a$ so that $|\xi_i|\leq a $   for
all $i$, then the  sum $S_n=\sum_{i=1}^{n}\xi_i$  obeys the following Bennett  inequality  (see also Bernstein \cite{B27}): for all $x>0,$
\begin{eqnarray}
 \mathbf{P}\Big(S_n \geq x  \Big)
&\leq& B_1(x, a, v):=\left(\frac{v^2}{x a + v^2} \right)^{\frac{x}{a} + \frac{ v^2}{a^2} }e^\frac{x}{a}  \\
&\leq& B_2(x, a, v):= \exp \left\{-\frac{x^2}{ 2(v^2+ x a/3  )} \right\},
\end{eqnarray}
where $v^2$ is the variance of $S_n.$

 Freedman have extended Bennett's  result to the case of
discrete-time supermartingales with bounded jumps.
Let $ (\xi _i,\mathcal{F}_i)_{i=1,...,n}$  be a sequence of  supermartingale differences. Denote by $S_k=\sum_{i=1}^{k}\xi_i$ and $\left\langle S\right\rangle_k=\sum_{i=1}^k\mathbf{E}(\xi _i^2|\mathcal{F} _{i-1})$. The well-known Freedman's inequality \cite{FR75}  for supermartingales states that: if $\xi_{i} \leq a$ for a positive constant $a$, then, for all $x, v > 0,$
\begin{eqnarray}
 \mathbf{P}\Big(S_k \geq x\ \mbox{and}\ \left\langle S\right\rangle_k\leq v^2\ \mbox{for some}\ k\in[1,n] \Big)
&\leq&    B_1(x,a, v)  \ \  \label{freedman}\\
 &\leq&  B_2(x, a, v) . \label{freedmanxf2}
\end{eqnarray}
In  particular, when $(\xi_i)_{i=1,...,n}$ are independent, the bounds (\ref{freedman}) and (\ref{freedmanxf2})  reduce to the bounds of Bennett \cite{B62}  and Bernstein  \cite{B27}, respectively. Moreover,  Freedman's inequality implies that  the bounds (\ref{freedman}) and (\ref{freedmanxf2}) not only hold for $ S_n $ but even hold  for
the maximum of partial sums $\max_{1\leq k \leq n } S_k.$

Replacing the quadratic characteristic $\left\langle S\right\rangle_k$ by $$\textrm{H}_k^y:= \sum_{i=1}^k\Big(\mathbf{E}(\xi_i^2 |\mathcal{F}_{i-1}) +\xi_i^2\textbf{1}_{\{|\xi_i|> y\}} \Big),$$
Dzhaparidze and van Zanten \cite{Dz01} have established a generalization of Freedman's inequality for martingales with unbounded differences:
 for all $x, y, v > 0,$
\begin{eqnarray}
 \mathbf{P}\Big(S_k \geq x\ \mbox{and}\ \textrm{H}_k^y \leq v^2\ \mbox{for some}\ k\in[1,n] \Big)  &\leq&    B_1(x,y, v)  . \label{vzineq}
\end{eqnarray}
In particular, if $|\xi_{i}| \leq a$ for all $i$, it holds $\textrm{H}_k^a = \left\langle S\right\rangle_k,$ and then the inequality of Dzhaparidze and van Zanten (\ref{vzineq}) reduces to Freeman's inequality (\ref{freedman}), implying  inequality (\ref{freedmanxf2}).

However, if $(\xi_{i})$ are not bounded from below, inequality (\ref{vzineq}) does not imply Freeman's inequality (\ref{freedman}) in general.
To fill this gap, we propose replacing the random variable $\textrm{H}_k^y$ in  inequality (\ref{vzineq}) by a smaller one  $\textrm{G}_k^{y},$ where
\begin{eqnarray}\label{GKY}
\textrm{G}_k^{y}   =\sum_{i=1}^k   \left( \frac{}{} \mathbf{E}(\xi_i^2\textbf{\textrm{1}}_{\{\xi_i \leq y\}} |\mathcal{F}_{i-1}) + \xi_i^2\textbf{\textrm{1}}_{\{\xi_i> y\}}\right).
\end{eqnarray}
 Our Theorem \ref{th23} states that, for all $x, y \geq 0$ and $v>0$,
\begin{eqnarray}
 \mathbf{P}\left( S_k  \geq x\ \mbox{and}\ \textrm{G}_k^{y} \leq v^2 \ \mbox{for some}\ k\in[1,n]\right)
&\leq & B_1\left(x,y,v \right).  \label{thd2die}
\end{eqnarray}
Since $\textrm{G}_k^{y}\leq \textrm{H}_k^y,$
inequality (\ref{thd2die}) implies the inequality of  Dzhaparidze and van Zanten (\ref{vzineq}).
Moreover, if $\xi_i \leq a$ for all $i\in[1,n]$ (may not be bounded from below), it holds $\textrm{G}_k^{a} = \left\langle S\right\rangle_k$ for all $k\in[1,n]$, and then  (\ref{thd2die}) also implies Freeman's inequality (\ref{freedman}). This  fills  the gap.

In Theorem \ref{th24}, we give a generalization of (\ref{vzineq}) to the supermartingales with non-square-integrable differences. Write
\[
\textrm{G}_n(\beta)  =\sum_{i=1}^n \left( \frac{}{} \mathbf{E}\big( |\xi_i|^\beta  |\mathcal{F}_{i-1} \big) + |\xi_i|^\beta \right)
\]
for a constant $\beta \in (1,2). $
Then,  for all $x, v>0$,
\begin{eqnarray}
 \mathbf{P}\left( \max_{1\leq k \leq n} S_k  \geq x\  \mbox{and}\ \textrm{G}_n(\beta) \leq v^\beta  \right)
 \leq  \exp\left\{-C(\beta) \left(\frac{x}{v} \right)^\frac{\beta}{\beta -1} \right\},
\end{eqnarray}
where $C(\beta)= \beta ^{\frac1{1-\beta}} \left(1- \beta^{-1} \right).$

Under the additional assumption of conditional symmetry,  Sason \cite{S12}  gave an improvements of Freedman's inequality (\ref{freedmanxf2}). Sason proved that  if  $(\xi _i,\mathcal{F}_i)_{i\geq 1}$ is a sequence of conditionally symmetric martingale differences with $|\xi_{i}| \leq a$ for a positive constant $a$, then, for all $x, v>0$,
\begin{eqnarray}
&& \mathbf{P}\left( S_k  \geq x\ \mbox{and}\   \left\langle S\right\rangle_k \leq v^2 \ \mbox{for some}\ k\in [1, n] \right) \ \ \ \ \ \nonumber \\
&&  \ \ \ \ \ \ \ \ \ \   \ \ \ \ \  \leq \ B_{0}(x, a, v):=  \exp \left\{-\lambda x   + \left(\frac{\cosh( \lambda a)-1}{ a^2}  \right) v^2 \right\}   \label{mainserds}
\end{eqnarray}
and $B_{0}(x, a, v) \leq B_{1}(x, a, v),$ where
\[
\lambda=   \frac{1}{a} \log\left( \sqrt{1+\frac{x^2 a^2}{v^4}}+ \frac{x a}{v^2}\right).
\] 

In the spirit  of Dzhaparide and van Zanten \cite{Dz01}, we establish the following  generalization of Sason's inequality (\ref{mainserds}). Define
\begin{eqnarray}\label{MKY}
\textrm{M}_k^{y}   =\sum_{i=1}^k  \Big( \mathbf{E}(\xi_i^2\textbf{1}_{\{|\xi_i| \leq y\}} |\mathcal{F}_{i-1}) + \xi_i^2\textbf{1}_{\{|\xi_i|> y\}} \Big).
\end{eqnarray}
Then, for all $x, v>0$ and all $y\geq 0$,
\begin{eqnarray}
 \mathbf{P}\left( S_k  \geq x\ \mbox{and}\ \textrm{M}_k^{y} \leq v^2 \ \mbox{for some}\ k\in [1, n] \right)
&\leq &  B_{0}(x, y, v), \label{mainerss}
\end{eqnarray}
where, by convention,
$$B_0 \left(x,0,v \right)= \lim_{y\rightarrow 0+} B_0 \left(x,y,v \right)=\exp \Big\{-\frac{x^2}{2\, v^2} \Big\}$$
applied when $y=0.$
If $|\xi_i| \leq a$ for all $i,$ then $\textrm{M}_k^{a} = \left\langle S\right\rangle_k$ for all $k$ and   (\ref{mainerss}) reduces to Sason's inequality (\ref{mainserds}). Notice that when $y=0,$ inequality (\ref{mainerss})
is known as De la Pe\~{n}a's inequality \cite{D99}. Hence, our bound establishes a connection between the inequalities
 of De la Pe\~{n}a and Sason.

The paper is organized as follows. We present our main results in Section \ref{sec2} and the application to self-normalized deviations in Section \ref{sf}, and devote to the proofs of the main results in Sections  \ref{sec5}
- \ref{sec7}.

\section{Main results}\label{sec2}

Assume that we are given a sequence of real-valued  supermartingale differences $(\xi _i,\mathcal{F}_i)_{i=0,...,n}$
 defined on some probability space $(\Omega ,\mathcal{F},\mathbf{P})$, where $\xi _0=0 $ and $\{\emptyset, \Omega\}=%
\mathcal{F}_0\subseteq ...\subseteq \mathcal{F}_n\subseteq \mathcal{F}$ are
increasing $\sigma$-fields. So we have $\mathbf{E}(\xi_{i}|\mathcal{F}_{i-1})\leq 0, \  i=1,...,n, $ by definition. Set
\begin{equation}  \label{matingal}
 S_k=\sum_{i=1}^k\xi _i,\quad k=1,...,n.
\end{equation}
Then $S=(S_k,\mathcal{F}_k)_{k=1,...,n}$ is a supermartingale.
Let $\left\langle S\right\rangle $   be  the quadratic characteristic and $[S]$ be the squared variation of the
supermartingale $S:$
\begin{equation}
 \left\langle S\right\rangle
_k=\sum_{i=1}^k\mathbf{E}(\xi _i^2|\mathcal{F} _{i-1})\ \ \ \ \ \textrm{and} \ \ \ \ \ \ \  [S]_k=\sum_{i=1}^{k}\xi_i^2.
\end{equation}

The  following theorem strengthens the inequality of Dzhaparide and van Zanten \cite{Dz01}.
\begin{theorem}
\label{th23} Assume $\mathbf{E} \xi _i^{2} < \infty$  for all $i\in [1,n]$.
Then, for all $x, y \geq 0$ and $v>0$,
\begin{eqnarray}
 \mathbf{P}\left( S_k  \geq x\ \mbox{and}\ \textrm{G}_k^{y} \leq v^2 \ \mbox{for some}\ k\in[1,n]\right)
&\leq & B_1\left(x,y,v \right)  \label{th2ie}\\
&\leq& B_2\left(x, y  ,v \right).\label{th2ie2}
\end{eqnarray}
where $\textrm{G}_k^{y} $ is defined by (\ref{GKY})  and $$B_1 \left(x,0,v \right)= \lim_{y\rightarrow 0+} B_1 \left(x,y,v \right)=B_2\left(x,0,v \right)$$
applied when $y=0$.
\end{theorem}

Since $\textrm{G}_k^{0} \leq   \left\langle S\right\rangle
_n + [S]_n$ for all $k\in[1,n],$ inequality (\ref{th2ie2}) implies the following result: for all $x, v> 0,$
\begin{eqnarray}\label{dsdsdf}
 \mathbf{P}\left( \max_{1\leq k \leq n} S_k  \geq x\ \mbox{and}\ \left\langle S\right\rangle
_n + [S]_n \leq v^2 \right)  \ \leq\ B_1\left(x,0,v \right).
\end{eqnarray}
This result slightly refines an earlier inequality of Bercu and Touati \cite{BT08}, where they   have obtained
the same bound on tail probabilities $\mathbf{P}\left(  S_n  \geq x\ \mbox{and}\ \left\langle S\right\rangle
_n + [S]_n \leq v^2 \right).$ Thus the sum $S_n $ has been strengthened to
the maximum of partial sums $\max_{1\leq k \leq n} S_k.$   A similar refinement is applied to   Delyon's inequality \cite{De09}, where he has established the following result, for all $x, v>0,$
\begin{eqnarray}\label{dfddf}
 \mathbf{P}\left( S_n  \geq x\ \mbox{and}\ \textrm{G}_n^{0} \leq v^2  \right)
&\leq& B_2\left(x, 0 ,v \right).
\end{eqnarray}

Consider the supermartingales with non-square-integrable  differences.
We have the following large deviation exponential bound, which can be regarded as a generalization of  Delyon's inequality  (\ref{dfddf}) or the inequality of Dzhaparide and van Zanten (\ref{vzineq}).

Denote by $x^+=\max\{x, 0\}$ and $x^-=-\min\{x, 0\}$  the positive and negative parts of $x,$ respectively.
\begin{theorem}
\label{th24} Assume $\mathbf{E} |\xi _i|^{\beta} < \infty$ for a constant $\beta \in (1,2)$ and  all $i\in [1,n]$. Write
\[
\textrm{G}^0_k(\beta)  =\sum_{i=1}^k \Big(  \mathbf{E}\big( (\xi_i^-)^\beta  |\mathcal{F}_{i-1} \big) + (\xi_i^+)^\beta \Big),\ \ \ \ k\in [1,n].
\]
Then, for all $x, v>0$,
\begin{eqnarray}
 \mathbf{P}\left( S_k  \geq x\ \mbox{and}\ \textrm{G}^0_k(\beta) \leq v^\beta \ \mbox{for some}\ k\in[1,n]\right)
 \leq   \exp\left\{-C(\beta) \left(\frac{x}{v} \right)^\frac{\beta}{\beta -1} \right\}, \label{f71}
\end{eqnarray}
where
\begin{eqnarray}\label{cbeta}
C(\beta)= \beta ^{\frac1{1-\beta}} \left(1- \beta^{-1} \right).
\end{eqnarray}
\end{theorem}

In particular, if $||\textrm{G}_n(\beta)||_{\infty}=O(n)$ as $n\rightarrow \infty$, then, for any $x>0,$
\begin{eqnarray}
 \mathbf{P}\left( \max_{1\leq k \leq n} S_k  \geq n x  \right)
 = O\bigg( \exp\Big\{- n \, C_x(\beta) \Big\} \bigg),
\end{eqnarray}
where $C_x(\beta) >0$ does not depend on $n.$

It is also interesting to see that when $\beta$  decreases  to $1$ in (\ref{f71}), the power $\frac{\beta}{\beta -1}$ is increasing to infinity and the corresponding constant $C(\beta)$ is decreasing  to $0.$  This means the larger the power, the smaller the  corresponding constant.

 One calls $(\xi _i,\mathcal{F}_i)_{i\geq 1}$  a sequence of \emph{conditionally symmetric} martingale differences, if $\mathbf{E}(\xi_i  >y |\mathcal{F}_{i-1})= \mathbf{E}(\xi_i < -y |\mathcal{F}_{i-1})$ for all $i$ and  any $y \geq0$.
Motivated by the result of Dzhaparide and van Zanten \cite{Dz01}, we  give a generalization of Sason's inequality (\ref{mainserds}) to the martingales with unbounded differences.
\begin{theorem}
\label{th3} Assume that $(\xi _i,\mathcal{F}_i)_{i\geq 1}$ is a sequence of conditionally symmetric martingale differences.
Then, for all $x, v>0$ and all $y\geq 0$,
\begin{eqnarray}
 \mathbf{P}\left( S_k  \geq x\ \mbox{and}\ \textrm{M}_k^{y} \leq v^2 \ \mbox{for some}\ k\in [1, n] \right)
&\leq &  B_{0}(x, y, v), \label{mainers}
\end{eqnarray}
where $\textrm{M}_k^{y}$ is defined by (\ref{MKY}).
\end{theorem}

If the martingale differences are bounded $|\xi_i|\leq a$ for a positive constant $a,$ then   $\textrm{M}_k^{a}= \left\langle S\right\rangle_k$ and  inequality (\ref{mainers}) with $y=a$
reduces to Sason's inequality (\ref{mainserds}). As pointed out by Sason \cite{S12}, inequality (\ref{mainers}) is the best possible that can be obtained from  Chernoff's inequality  $\mathbf{P}\left(  S_n \geq x \right) \leq  \inf_{\lambda\geq 0} \mathbf{E}e^{\lambda(S_n- x)}$
 under the present assumption in a
 certain sense. Indeed, if  $(\xi_i)_{i \geq1}$ are i.i.d.\ random variables and satisfy the following  distribution
\begin{eqnarray}
\mathbf{P}(\xi_i=y  )= \mathbf{P}(\xi_i= -y )= \frac{v^2}{2ny^2}  \ \ \ \ \mbox{and} \ \ \ \mathbf{P}(\xi_i= 0 )= 1- \frac{v^2}{ ny^2},
\end{eqnarray}
then the bound (\ref{mainers}) equals to $\lim_{n\rightarrow \infty}\inf_{\lambda\geq 0} \mathbf{E}e^{\lambda(S_n- x)}$.

Since  $\lim_{y\rightarrow 0}\lambda = \frac{x}{v^2}$ and $\textrm{G}_k^{0}=[S]_k, $ inequality (\ref{mainers}) reduces to De la Pe\~{n}a's  inequality \cite{D99}  as $y\rightarrow 0$: for all $x,  v>0$,
\begin{eqnarray}\label{equdela}
 \mathbf{P}\left( S_k  \geq x\ \mbox{and}\ [S]_k \leq v^2 \ \mbox{for some}\ k\in \mathbb{N}\right)
&\leq &   B_{0}(x, 0, v).
\end{eqnarray}
 Thus inequality (\ref{mainers}) establishes the connection between the inequalities of  De la Pe\~{n}a and Sason.

\section{Application to self-normalized deviations}\label{sf}
As an application of Theorem \ref{th24},  consider the self-normalized deviations for independent random variables.
\begin{theorem}\label{co01}
Assume  that $(\xi _i)_{i=1,...,n}$ is a sequence of independent and symmetric random variables. Denote by $$\textrm{V}_n(\beta) =\Big(\sum_{i=1}^{n} |\xi_i|^\beta \Big)^{1/\beta}$$
for a constant $\beta \in (1,2].$  Then, for all $x > 0$,
\begin{eqnarray}
 \mathbf{P}\left( \max_{1\leq k \leq n} S_k / \textrm{V}_n(\beta) \geq x\right)
 \leq   \exp\bigg\{- \widetilde{C} (\beta) x^\frac{\beta}{\beta -1}   \bigg\}, \label{fdsdgsh}
\end{eqnarray}
where
\begin{eqnarray}
\widetilde{C} (\beta)= \Big(\frac \beta 2 \Big) ^{\frac1{1-\beta}} \left(1- \beta^{-1} \right).
\end{eqnarray}
\end{theorem}
\noindent\emph{Proof of Theorem \ref{co01}.}
Assume that $(\xi _i)_{i=1,...,n}$ are  independent and  symmetric. Set $$\mathcal{F}_{i}=\sigma \left\{\xi_{k}, k\leq i, |\xi_j|, 1\leq j \leq n \right\}.$$ Since $\xi _i$ is symmetric, it is easy to see that
\[
\mathbf{E}\Big( \frac{\xi_i}{\textrm{V}_n(\beta)} \Big|\mathcal{F}_{i-1} \Big)=\mathbf{E}\Big(\xi_i \Big| |\xi_{i}| \Big) \frac{1}{\textrm{V}_n(\beta)}  =0.
\]
Therefore, $\big(\xi_i /  \textrm{V}_n(\beta),\mathcal{F}_i\big)_{i=1,...,n}$ is a sequence of   martingale differences.
Notice that
\begin{eqnarray*}
\frac{1}{\textrm{V}_n(\beta)^\beta}\sum_{i=1}^k \Big(  \mathbf{E}\big( (\xi_i^-)^\beta  |\mathcal{F}_{i-1} \big) + (\xi_i^+)^\beta \Big) &\leq & \frac{1}{\textrm{V}_n(\beta)^\beta}\sum_{i=1}^n \Big(  \mathbf{E}\big( |\xi_i|^\beta  |\mathcal{F}_{i-1} \big) + |\xi_i|^\beta \Big) \\
&=& \frac{2}{\textrm{V}_n(\beta)^\beta}\sum_{i=1}^{n}|\xi_i|^\beta     = 2.
\end{eqnarray*}
Applying Theorem \ref{th24}   to $\big(\xi_i /  \textrm{V}_n(\beta),\mathcal{F}_i\big)_{i=1,...,n}$, we obtain (\ref{fdsdgsh}).\hfill\qed

The power $x^\frac{\beta}{\beta -1}$ in (\ref{fdsdgsh})  is the best possible
for $x$ in the moderate deviation and large deviation ranges. Indeed, Jing, Liang  and Zhou
\cite{JLZ12}  have obtained the following self-normalized moderate deviation result. Assume that  $$\mathbf{P}(\xi_i \geq x) = \mathbf{P}(\xi_i \leq - x) \sim \frac{c}{x^\alpha} h_i(x), \ \ \ \ \ x \rightarrow \infty, $$
where $\alpha \in (0, 2),$ $c >0$ and $h_i(x)$'s are slowly varying at $\infty.$   Under certain conditions on the tail probabilities of $\xi_i$ (cf.\ Theorem 2.3 of \cite{JLZ12} for details), for $x_n \rightarrow \infty$ and $x_n= o(n^{(\beta-1)/\beta})$ and $\beta > \max\{1,  \alpha\},$  the limit exists
\begin{eqnarray}\label{sdbof}
 \lim_{n\rightarrow \infty} x_{n}^{-\frac{\beta}{\beta -1}} \log    \mathbf{P}\left(   S_n /  \textrm{V}_n(\beta)  \geq x_n \right) = - (\beta -1)C_\alpha(\beta),
\end{eqnarray}
where $C_\alpha(\beta)$ is a positive constant depending on  $\alpha$ and $\beta.$ Equality (\ref{sdbof}) suggests that the power $x^\frac{\beta}{\beta -1}$ in (\ref{fdsdgsh})  is the best possible for moderate $x$'s.
See also Shao \cite{S97} for self-normalized large deviation result.

\section{Proof of Theorem \ref{th23}} \label{sec5}
Assume $(\xi _i,\mathcal{F}_i)_{i=0,...,n}$ a sequence of square integrable supermartingale differences. For any
nonnegative numbers $y$ and $\lambda$, define the exponential multiplicative
martingale $Z(\lambda )=(Z_k(\lambda ),\mathcal{F}_k)_{k=0,...,n},$ where
\[
Z_k(\lambda )=\prod_{i=1}^k\frac{\exp\left\{\lambda \xi_i -\frac{1}{2} (\lambda\xi_i)^2  \textbf{1}_{\{\xi_i> y\}}   \right\}}{\mathbf{E}\left(\exp\left\{\lambda \xi_i -\frac{1}{2} (\lambda\xi_i)^2  \textbf{1}_{\{\xi_i> y\}}   \right\} |
\mathcal{F}_{i-1} \right)},\quad \quad  \quad Z_0(\lambda )=1.
\]
If $T$ is a stopping time, then $Z_{T\wedge k}(\lambda )$, $\lambda >0$, is also a martingale, where
\[
Z_{T\wedge k}(\lambda )=\prod_{i=1}^{T\wedge k}\frac{\exp\left\{\lambda \xi_i -\frac{1}{2} (\lambda\xi_i)^2  \textbf{1}_{\{\xi_i> y\}}   \right\}}{\mathbf{E}\left(\exp\left\{\lambda \xi_i -\frac{1}{2} (\lambda\xi_i)^2  \textbf{1}_{\{\xi_i> y\}}   \right\} |
\mathcal{F}_{i-1} \right)}, \quad \quad Z_0(\lambda )=1. \label{C-1}
\]
Then for any nonnegative number $\lambda,$ we have the
following conjugate probability measure $\mathbf{P}_\lambda $ on $(\Omega ,\mathcal{F})$:
\begin{equation}
d\mathbf{P}_\lambda =Z_{T\wedge n}(\lambda )d\mathbf{P}.  \label{chmeasure1}
\end{equation}

\begin{lemma}
\label{lemma2} For all $y\geq 0$ and  all $\lambda>0,$ it holds
\begin{eqnarray*}
 \mathbf{E}\left( \exp\left\{\lambda \xi_i -\frac{1}{2} (\lambda\xi_i)^2  \textbf{\emph{1}}_{\{\xi_i> y\}}   \right\}\Bigg|\mathcal{F}_{i-1} \right)  \leq  \exp\Bigg\{ \left( \frac{e^{\lambda y}-1-\lambda y }{y^2} \right)\mathbf{E}(\xi_i^2\textbf{\emph{1}}_{\{\xi_i\leq y\}}  |\mathcal{F}_{i-1} ) \Bigg\},
\end{eqnarray*}
where, by convention, $\frac{e^{\lambda y}-1-\lambda y }{y^2}=\frac{\lambda^2}{2}$ applied when $y=0$.
\end{lemma}
\noindent\emph{Proof.} Let $y\geq0.$
If $\xi_i\leq y$, since the function
$$g(x)=\frac{e^x-1-x}{x^2}$$
is increasing in $x \in \mathbf{R}$ (by convention $g(0)=1/2$),  we have, for all $\lambda>0$,
\begin{eqnarray}\label{hjk125}
\frac{e^{\lambda \xi_i}-1- \lambda \xi_i}{ (\lambda \xi_i)^2}  \leq   \frac{e^{\lambda y}-1- \lambda y}{ (\lambda y)^2}.
\end{eqnarray}
If $\xi_i > y$, since  $\exp\left\{x -\frac{1}{2} x^2   \right\}  \leq  1 + x$ for all $x\geq 0$, it follows that, for all $\lambda>0$,
\begin{eqnarray}\label{hjk226}
 \exp\left\{\lambda \xi_i -\frac{1}{2} (\lambda\xi_i)^2  \right\}  \leq   1 + \lambda \xi_i .
\end{eqnarray}
Combining (\ref{hjk125}) and (\ref{hjk226}) together, we find that, for all  $y\geq 0$ and all $\lambda>0,$
\begin{eqnarray*}
 \exp\left\{\lambda \xi_i -\frac{1}{2} (\lambda\xi_i)^2  \textbf{1}_{\{\xi_i> y\}}   \right\}  \leq   1 + \lambda \xi_i+ \left( \frac{e^{\lambda y}-1-\lambda y }{y^2} \right)\xi_i^2\textbf{1}_{\{\xi_i\leq y\}}.
\end{eqnarray*}
Taking conditional expectations on both sides of the last inequality, we deduce that
\begin{eqnarray}\label{sfsdsv}
 \mathbf{E}\left( \exp\left\{\lambda \xi_i -\frac{1}{2} (\lambda\xi_i)^2  \textbf{1}_{\{\xi_i> y\}}   \right\}\Bigg|\mathcal{F}_{i-1} \right)  \leq   1 + \left( \frac{e^{\lambda y}-1-\lambda y }{y^2} \right)\mathbf{E}(\xi_i^2\textbf{1}_{\{\xi_i\leq y\}}  |\mathcal{F}_{i-1} ).
\end{eqnarray}
Using the inequality $1+x \leq e^x$ in the right-hand side of (\ref{sfsdsv}),
we obtain the desired inequality.\hfill\qed

\vspace{0.3cm}

\noindent\emph{Proof of Theorem \ref{th23}.}
For any $x,  v>0$ and any $y\geq0$, define the stopping time $T$:
\[
T(x,y,v) =\min\{k\in [1, n]: S_k\geq x\ \mbox{and}\ G_k^{y}\leq v^2 \},
\]
with the convention that $\min{\emptyset}=0$. Then it follows that
 \[
\textbf{1}_{\{ S_k\geq x\  \mbox{and}\ G_k^{y}\leq v^2\ \mbox{for some}\ k\in[1,n] \}} = \sum_{k=1}^{n}  \textbf{1}_{\{ T(x,y,v) =k\}}.
\]
Denote by $\mathbf{E}_{\lambda}$ the expectation with
respect to the conjugate probability measure  $\mathbf{P}_\lambda$.
Using the change of probability measure (\ref{chmeasure1}), we have, for all $x,\lambda, v>0$ and all $y\geq 0$,
\begin{eqnarray}
  && \mathbf{P}( S_k\geq x\ \mbox{and}\  G_k^{y}\leq v^2\ \mbox{for some}\ k\in[1,n]) \nonumber\\
 &=& \mathbf{E}_{\lambda} \left( Z_{T\wedge n}(\lambda)^{-1}\textbf{1}_{\{S_k\geq x\  \mbox{and}\  G_k^{y}\leq v^2\ \mbox{for some}\ k\in[1,n]\}} \right) \nonumber \\
 &=& \sum_{k=1}^{n}\mathbf{E}_{\lambda} \left(\exp\left\{-\lambda S_{k}+\frac{\lambda^2}{2} [S]_{k} (y)+ \widetilde{\Psi}_{k}(\lambda) \right\} \textbf{1}_{\{T(x,y,v) =k\}} \right), \label{ghna}
\end{eqnarray}
where
$$[S]_{k} (y)=\sum_{i=1}^{k}\xi_i^2\textbf{1}_{\{\xi_i> y \}} $$
and
\[
\widetilde{\Psi}_{k}(\lambda)= \sum_{i=1}^k \log \mathbf{E} \left(\exp\left\{\lambda \xi_i -\frac{1}{2} (\lambda\xi_i)^2  \textbf{1}_{\{\xi_i> y\}}   \right\} \bigg|\mathcal{F}_{i-1} \right).
\]
Since the function $g(x)$ in increasing in $x$ and $g(0)=1/2,$ we have $$\frac{\lambda^2}{2} \leq   \frac{e^{\lambda y}-1- \lambda y}{ y^2}\  \ \ \ \ \textrm{for all}\ \  y, \lambda>0 .$$
Hence, by Lemma \ref{lemma2} and the last inequality,
\begin{eqnarray}
  && \mathbf{P}( S_k\geq x\ \mbox{and}\  G_k^{y} \leq v^2\ \mbox{for some}\ k\in[1,n]) \nonumber\\
 &\leq& \sum_{k=1}^{n}\mathbf{E}_{\lambda} \left(\exp \left\{-\lambda S_{k} +\frac{\lambda^2}{2} [S]_{k} (y)  + \left(\frac{e^{\lambda y}-1- \lambda y}{ y^2}  \right) \langle S\rangle_{k} (y) \right\} \textbf{1}_{\{T(x,y,v) =k\}} \right)\nonumber\\
 &\leq& \sum_{k=1}^{n}\mathbf{E}_{\lambda} \left(\exp \left\{-\lambda S_{k}   + \left(\frac{e^{\lambda y}-1- \lambda y}{ y^2}  \right) G_k^{y}\right\} \textbf{1}_{\{T(x,y,v) =k\}} \right), \label{dgcfj}
\end{eqnarray}
where  $\langle S\rangle_{k} (y)=\sum_{i=1}^k\mathbf{E}(\xi _i^2\textbf{1}_{\{ \xi_i \leq y \}}|\mathcal{F} _{i-1})$.
Therefore, by the fact $S_{k}\geq x$ and $G_k^{y} \leq v^2$  on the set $\{T(x,y,v) =k\}$, inequality (\ref{dgcfj}) implies that,
for all $x,\lambda, v>0$ and all $y\geq 0$,
\begin{eqnarray}
  && \mathbf{P}( S_k\geq x\ \mbox{and}\  G_k^{y} \leq v^2\ \mbox{for some}\ k\in[1,n]) \nonumber\\
 &\leq& \sum_{k=1}^{n}\mathbf{E}_{\lambda} \left(\exp \left\{-\lambda x   + \left(\frac{e^{\lambda y}-1- \lambda y}{ y^2}  \right) v^2 \right\} \textbf{1}_{\{T(x,y,v) =k\}} \right)\nonumber\\
 &\leq&  \exp \left\{-\lambda x   + \left(\frac{e^{\lambda y}-1- \lambda y}{ y^2}  \right) v^2 \right\}  \label{llastie1}.
\end{eqnarray}
The last inequality attains its minimum at
\[
\lambda=  \lambda(x)= \frac{1}{y} \log\left( 1+ \frac{x y}{v^2}\right).
\]
Substituting $\lambda = \lambda(x)$ in (\ref{llastie1}), we obtain (\ref{th2ie}). Using the inequality
 \[  e^{t}-1- t \leq \frac{t^2}{2(1-\frac t 3)},\ \ \  \ \ t\geq 0, \]
 we get, for all $x, v>0$ and all $y\geq 0$,
\begin{eqnarray}
 \inf_{ \lambda> 0 } \exp \left\{-\lambda x   + \left(\frac{e^{\lambda y}-1- \lambda y}{ y^2}  \right) v^2 \right\} &\leq&  \inf_{ \lambda> 0 } \exp \left\{-\lambda x   + \frac{\lambda^2v^2}{2(1-\frac{\lambda y} 3)}   \right\}\nonumber\\
 &\leq& B_2\left(x,y,v \right). \nonumber
\end{eqnarray}
Thus  (\ref{th2ie}) implies (\ref{th2ie2}).
This completes the proof  of Theorem \ref{th23}.\hfill\qed

\section{Proof of Theorem \ref{th24}} \label{sec6}
Assume $\mathbf{E}|\xi_i|^\beta < \infty$ for a constant $\beta \in (1,2)$ and for all $i \in [1, n]$. For any nonnegative numbers $\lambda$, define the exponential multiplicative
martingale $Z(\lambda )=(Z_k(\lambda ),\mathcal{F}_k)_{k=0,...,n},$ where
\[
Z_k(\lambda )=\prod_{i=1}^k\frac{\exp\left\{\lambda \xi_i -  (\lambda \xi_i^+)^\beta     \right\}}{\mathbf{E}\left(\exp\left\{\lambda \xi_i -(\lambda \xi_i^+)^\beta   \right\} |
\mathcal{F}_{i-1} \right)},\quad \quad  \quad Z_0(\lambda )=1.
\]
If $T$ is a stopping time, then $Z_{T\wedge k}(\lambda )$, $\lambda \geq 0$, is also a martingale, where
\[
Z_{T\wedge k}(\lambda )=\prod_{i=1}^{T\wedge k}\frac{\exp\left\{\lambda \xi_i -  (\lambda \xi_i^+)^\beta    \right\}}{\mathbf{E}\left(\exp\left\{\lambda \xi_i -  (\lambda \xi_i^+)^\beta   \right\} |
\mathcal{F}_{i-1} \right)}, \quad \quad Z_0(\lambda )=1.
\]
Then for any nonnegative number $\lambda,$ we introduce the
following conjugate probability measure $\mathbf{P}_\lambda $ on $(\Omega ,\mathcal{F})$:
\begin{equation}
d\mathbf{P}_\lambda =Z_{T\wedge n}(\lambda )d\mathbf{P}.  \label{chmeas}
\end{equation}
\begin{lemma}
\label{lemmas2t} If  $\mathbf{E}|\xi_i|^\beta < \infty$ for a constant $\beta \in (1,2)$, then, for all $\lambda>0,$
\begin{eqnarray*}
 \mathbf{E}\Big( \exp\left\{\lambda \xi_i - \lambda^\beta ( \xi_i^+)^\beta  \right\} \Big|\mathcal{F}_{i-1} \Big)  \leq  \exp\bigg\{  \lambda^\beta \mathbf{E}( ( \xi_i^-)^\beta \big|\mathcal{F}_{i-1})\bigg\}.
\end{eqnarray*}
\end{lemma}
\noindent\emph{Proof.} It is easy to see that, for all $x \in \mathbf{R}$ and $\beta \in (1,2)$,
\[
\exp\left\{x-(x^+)^\beta \right\} \leq 1+ x +(x^-)^\beta .
\]
With $x=\lambda \xi_i$, we easily obtain, for all $\lambda\geq 0$,
\begin{eqnarray}
\exp\left\{\lambda \xi_i -  (\lambda \xi_i^+)^\beta  \right\}  \leq   1 + \lambda \xi_i +  (\lambda \xi_i^-)^\beta . \label{phat}
\end{eqnarray}
Taking conditional expectations on both sides of (\ref{phat}), we get
\begin{eqnarray*}
 \mathbf{E}\Big( \exp\left\{\lambda \xi_i - \lambda^\beta ( \xi_i^+)^\beta  \right\} \Big|\mathcal{F}_{i-1} \Big)  \leq  1 +  \lambda^\beta \mathbf{E}( ( \xi_i^-)^\beta \big|\mathcal{F}_{i-1}) .
\end{eqnarray*}
Using the inequality $ 1+x \leq e^x$, we obtain the desired inequality.  \hfill\qed

\vspace{0.3cm}

\noindent \emph{Proof of Theorem \ref{th24}.}
For given $x,v>0$, define the stopping time $T:$
\[
T=\min\{k\in [1, n]: S_k\geq x\ \mbox{and}\ \textrm{G}^0_k(\beta) \leq v^\beta \},
\]
with the convention that $\min{\emptyset}=0$. Then we have
 \[
\textbf{1}_{\{ S_k\geq x\  \mbox{and}\   G^0_{k} (\beta) \leq v^\beta\ \mbox{for some}\ k\in[1,n] \}} = \sum_{k=1}^{n}  \textbf{1}_{\{ T =k\}}.
\]
Denote by $\mathbf{E}_{\lambda}$ the expectation with
respect to the conjugate probability measure (\ref{chmeas}).
Using the change of measure (\ref{chmeas}), we get, for all $x,\lambda, v>0$,
\begin{eqnarray}
  && \mathbf{P}( S_k\geq x\ \mbox{and}\  \textrm{G}^0_k(\beta)\leq v^\beta\ \mbox{for some}\ k\in[1,n]) \nonumber\\
 &=& \mathbf{E}_{\lambda} \Big( Z_{T\wedge n}(\lambda)^{-1}\textbf{1}_{\{S_k\geq x\  \mbox{and}\  G^0_k(\beta)\leq v^2\ \mbox{for some}\ k\in[1,n]\}} \Big) \nonumber \\
 &=& \sum_{k=1}^{n}\mathbf{E}_{\lambda} \Big( \exp\{-\lambda S_{k}+ \lambda^\beta \sum_{i=1}^k ( \xi_i^+)^\beta+ \widehat{\Psi}_{k}(\lambda) \} \textbf{1}_{\{T =k\}} \Big), \label{ghnsa}
\end{eqnarray}
where
\[
\widehat{\Psi}_{k}(\lambda)= \sum_{i=1}^k \log \mathbf{E} \exp\left\{\lambda \xi_i - (\lambda\xi_i^+)^\beta    \right\}.
\]
From inequality (\ref{ghnsa}), by Lemma \ref{lemmas2t}, it follows that, for all $x,\lambda, v>0$,
\begin{eqnarray}
 && \mathbf{P}( S_k\geq x\ \mbox{and}\  \textrm{G}^0_k(\beta)\leq v^\beta\ \mbox{for some}\ k\in[1,n])  \ \ \ \ \ \ \   \nonumber\\
 && \ \ \ \ \ \ \  \ \ \ \ \ \ \  \leq \ \sum_{k=1}^{n}\mathbf{E}_{\lambda} \Big( \exp\{-\lambda S_{k}+ \lambda^\beta \textrm{G}^0_k(\beta) \} \textbf{1}_{\{T =k\}} \Big). \nonumber
\end{eqnarray}
Since $S_k\geq x$ and $\textrm{G}^0_k(\beta)\leq v^\beta$ on the set $\{T =k\}$, we obtain, for all $x,\lambda, v>0$,
\begin{eqnarray}
 \mathbf{P}( S_k\geq x\ \mbox{and}\  \textrm{G}^0_k(\beta)\leq v^\beta\ \mbox{for some}\ k\in[1,n])
  \leq    \exp\{-\lambda x+ \lambda^\beta v^\beta \} . \label{sfgdf}
\end{eqnarray}
The last inequality attains its minimum at
\[
\lambda=  \lambda(x)= \left( \frac{x}{ \beta v^\beta}\right)^{\frac{1}{ \beta -1}} .
\]
Substituting $\lambda = \lambda(x)$ in (\ref{sfgdf}), we get the desired inequality.\hfill\qed

\section{Proof of Theorem \ref{th3}} \label{sec7}
Assume that $(\xi _i,\mathcal{F}_i)_{i=0,...,n}$ is a sequence of  martingale differences. For any
nonnegative numbers $y$ and $\lambda$, define the exponential multiplicative
martingale $Z(\lambda )=(Z_k(\lambda ),\mathcal{F}_k)_{k=0,...,n},$ where
\[
Z_k(\lambda )=\prod_{i=1}^k\frac{\exp\left\{\lambda \xi_i -\frac{1}{2} (\lambda\xi_i)^2  \textbf{1}_{\{|\xi_i|> y\}}   \right\}}{\mathbf{E}\left(\exp\left\{\lambda \xi_i -\frac{1}{2} (\lambda\xi_i)^2  \textbf{1}_{\{|\xi_i|> y\}}   \right\} |
\mathcal{F}_{i-1} \right)},\quad \quad  \quad Z_0(\lambda )=1.
\]
If $T$ is a stopping time, then $Z_{T\wedge k}(\lambda )$, $\lambda >0$, is also a martingale, where
\[
Z_{T\wedge k}(\lambda )=\prod_{i=1}^{T\wedge k}\frac{\exp\left\{\lambda \xi_i -\frac{1}{2} (\lambda\xi_i)^2  \textbf{1}_{\{|\xi_i|> y\}}   \right\}}{\mathbf{E}\left(\exp\left\{\lambda \xi_i -\frac{1}{2} (\lambda\xi_i)^2  \textbf{1}_{\{|\xi_i|> y\}}   \right\} |
\mathcal{F}_{i-1} \right)}, \quad \quad Z_0(\lambda )=1. \label{C-1}
\]
Then for any nonnegative number $\lambda,$ we have the
following conjugate probability measure $\mathbf{P}_\lambda $ on $(\Omega ,\mathcal{F})$:
\begin{equation}
d\mathbf{P}_\lambda =Z_{T\wedge n}(\lambda )d\mathbf{P}.  \label{chmeasure2}
\end{equation}

\begin{lemma}  Assume that $(\xi _i,\mathcal{F}_i)_{i=0,...,n}$ is a sequence of conditionally symmetric martingale differences. For all $\lambda, y\geq 0,$   it holds
\begin{eqnarray*}
 \mathbf{E}\left( \exp\left\{\lambda \xi_i -\frac{1}{2} (\lambda\xi_i)^2  \textbf{\emph{1}}_{\{|\xi_i|> y\}}   \right\}\Bigg|\mathcal{F}_{i-1} \right)  \leq  \exp\bigg\{\left( \frac{\cosh(\lambda y)-1}{y^2} \right)\mathbf{E}(\xi_i^2\textbf{\emph{1}}_{\{|\xi_i|\leq y\}}  |\mathcal{F}_{i-1} )\bigg\},
\end{eqnarray*}
where by convention $\frac{\cosh(\lambda y)-1}{y^2}=\frac{\lambda^2}{2}$ when $y=0$.
\end{lemma}
\noindent\emph{Proof.} Let $y\geq0.$
When $|\xi_i|\leq y$,     it follows that $  \xi_i^{2k}   \leq y^{2k-2}  \xi_i^{2}  $ and that
\begin{eqnarray}
\cosh(\lambda \xi_i)     \leq   1+  \frac{  \xi_i^{2}  }{y^2} \sum_{k=1}^{\infty}\frac{ (\lambda y)^{2k}}{(2k)!}
   =   1+  \frac{  \xi_i^{2}  }{y^2} \Big( \cosh(\lambda y)-1 \Big).\label{hjk1}
\end{eqnarray}
When $|\xi_i| > y$, since  $\cosh(x) \leq   \exp\left\{\frac{1}{2} x^2   \right\} $ for all $x\in \mathbf{R}$, it follows that, for all $\lambda>0$,
\begin{eqnarray}\label{hjk2}
 \Big( \cosh(\lambda \xi_i)  \Big) \exp\left\{ -\frac{1}{2} (\lambda\xi_i )^2  \right\}  \leq   1   .
\end{eqnarray}
Combining (\ref{hjk1}) and (\ref{hjk2}) together, we find that, for all  $\lambda, y\geq 0,$
\begin{eqnarray*}
 \Big( \cosh(\lambda \xi_i) \Big) \exp\left\{  -\frac{1}{2} (\lambda\xi_i)^2  \textbf{1}_{\{|\xi_i|> y\}}   \right\}  \leq   1 + \lambda \xi_i+ \left( \frac{e^{\lambda y}-1-\lambda y }{y^2} \right)\xi_i^2\textbf{1}_{\{|\xi_i|\leq y\}}.
\end{eqnarray*}
Taking conditional expectations on both sides of the last inequality, we have, for all $\lambda, y\geq 0,$
\begin{eqnarray}
 && \mathbf{E}\left( \Big( \cosh(\lambda \xi_i) \Big) \exp\left\{  -\frac{1}{2} (\lambda\xi_i)^2  \textbf{1}_{\{|\xi_i|> y\}}   \right\}\Bigg|\mathcal{F}_{i-1} \right) \nonumber \\
  &&\ \ \ \ \ \ \ \ \ \ \ \ \ \ \ \ \ \ \ \  \ \ \ \ \ \ \ \ \ \ \ \ \ \   \leq   1 + \left( \frac{\cosh(\lambda y)-1}{y^2} \right)\mathbf{E}(\xi_i^2\textbf{1}_{\{|\xi_i|\leq y\}}  |\mathcal{F}_{i-1} ). \label{fine9}
\end{eqnarray}
Since $(\xi _i,\mathcal{F}_i)_{i=0,...,n}$ are conditionally symmetric, it holds
\begin{eqnarray}
  & &\mathbf{E}\left( \Big( \exp\{\lambda \xi_i\} \Big) \exp\left\{  -\frac{1}{2} (\lambda\xi_i)^2  \textbf{1}_{\{|\xi_i|> y\}}   \right\}\Bigg|\mathcal{F}_{i-1} \right) \nonumber \\  & & \ \ \ \ \ \ \ \ \   \ \ \   \   \ \   \ \ \ \  \ \ \ \ \ =   \mathbf{E}\left( \Big( \exp\{-\lambda \xi_i\} \Big) \exp\left\{  -\frac{1}{2} (\lambda\xi_i)^2  \textbf{1}_{\{|\xi_i|> y\}}   \right\}\Bigg|\mathcal{F}_{i-1} \right). \label{sfcgh}
\end{eqnarray}
Note that $\cosh(\lambda \xi_i)=\frac12 (\exp\{\lambda \xi_i\}+ \exp\{-\lambda \xi_i\} ).$
Hence (\ref{sfcgh}) implies that
\begin{eqnarray}
  & & \mathbf{E}\left( \Big( \cosh(\lambda \xi_i) \Big) \exp\left\{  -\frac{1}{2} (\lambda\xi_i)^2  \textbf{1}_{\{|\xi_i|> y\}}   \right\}\Bigg|\mathcal{F}_{i-1} \right)\nonumber \\  & &\ \ \ \ \ \ \ \ \ \   \ \ \ \ \ \ \ \   \ \   \ \ \ \  \ \ \ \ \ \   =\ \mathbf{E}\left( \Big( \exp\{\lambda \xi_i\} \Big) \exp\left\{  -\frac{1}{2} (\lambda\xi_i)^2  \textbf{1}_{\{|\xi_i|> y\}}   \right\}\Bigg|\mathcal{F}_{i-1} \right). \label{dfd1s}
\end{eqnarray}
Combining (\ref{fine9}) and (\ref{dfd1s}) together, we obtain
\begin{eqnarray*}
 \mathbf{E}\left( \exp\left\{\lambda \xi_i -\frac{1}{2} (\lambda\xi_i)^2  \textbf{1}_{\{|\xi_i|> y\}}   \right\}\Bigg|\mathcal{F}_{i-1} \right)  \leq   1 + \left( \frac{\cosh(\lambda y)-1}{y^2} \right)\mathbf{E}(\xi_i^2\textbf{ 1}_{\{|\xi_i|\leq y\}}  |\mathcal{F}_{i-1} ).
\end{eqnarray*}
Using the inequality $1+x \leq e^x,$ we obtain
the  desired inequality.
\hfill\qed

\vspace{0.3cm}

\noindent\emph{Proof of Theorem \ref{th3}.}
For any $y\geq0$ and any $x, v>0$, define the stopping time $T$:
\[
T(x,y,v) =\min\{k\in [1, n]: S_k\geq x\ \mbox{and}\ G_k^{y}\leq v^2 \},
\]
with the convention that $\min{\emptyset}=0$. Then it follows that
 \[
\textbf{1}_{\{ S_k\geq x\  \mbox{and}\ G_k^{y}\leq v^2\ \mbox{for some}\ k\in[1,n] \}} = \sum_{k=1}^{n}  \textbf{1}_{\{ T(x,y,v) =k\}}.
\]
Denote $\mathbf{E}_{\lambda}$ the expectation with
respect to $\mathbf{P}_{\lambda}$.
Using the change of probability measure (\ref{chmeasure2}), we have, for all $x,\lambda, v>0$ and all $y\geq 0$,
\begin{eqnarray}
  && \mathbf{P}( S_k\geq x\ \mbox{and}\  G_k^{y}\leq v^2\ \mbox{for some}\ k\in[1,n]) \nonumber\\
 &=& \mathbf{E}_{\lambda} \left( Z_{T\wedge n}(\lambda)^{-1}\textbf{1}_{\{S_k\geq x\  \mbox{and}\  G_k^{y}\leq v^2\ \mbox{for some}\ k\in[1,n]\}} \right) \nonumber \\
 &=& \sum_{k=1}^{n}\mathbf{E}_{\lambda} \left(\exp\left\{-\lambda S_{k}+\frac{\lambda^2}{2} [S]_{k} (y)+   \breve{\Psi}  _{k}(\lambda) \right\} \textbf{1}_{\{T(x,y,v) =k\}} \right), \label{ghna}
\end{eqnarray}
where $[S]_{k} (y)=\sum_{i=1}^{k}\xi_i^2\textbf{1}_{\{|\xi_i|> y \}}$ and
\[
 \breve{\Psi} _{k}(\lambda)= \sum_{i=1}^k \log \mathbf{E} \left(\exp\left\{\lambda \xi_i -\frac{1}{2} (\lambda\xi_i)^2  \textbf{1}_{\{|\xi_i|> y \}}   \right\} \bigg|\mathcal{F}_{i-1} \right).
\]
Since the function $f(y)= \frac{\cosh( \lambda y)-1}{ y^2}$ is increasing in $y\geq 0$ and $f(0)=\lambda^2/2,$ we have $$\frac{\lambda^2}{2} \leq   \frac{\cosh( \lambda y)-1}{ y^2}\   \ \ \textrm{for all}\ \  y, \lambda>0 .$$
By Lemma \ref{lemma2},  we find that $$\widetilde{\Psi}_{k}(\lambda) \leq \left(\frac{\cosh( \lambda y)-1}{ y^2}  \right) \langle S\rangle_{k} (y),$$ where $\langle S\rangle_{k} (y)=\sum_{i=1}^k\mathbf{E}(\xi _i^2\textbf{1}_{\{|\xi_i|> y \}}|\mathcal{F} _{i-1})$.
Then, from  equality (\ref{ghna}), it follows that,  for all $x,\lambda, v>0$ and all $y\geq 0$,
\begin{eqnarray}
  && \mathbf{P}( S_k\geq x\ \mbox{and}\  G_k^{y} \leq v^2\ \mbox{for some}\ k\in[1,n]) \nonumber\\
 &\leq& \sum_{k=1}^{n}\mathbf{E}_{\lambda} \left(\exp \left\{-\lambda S_k +\frac{\lambda^2}{2} [S]_{k} (y)  + \left(\frac{\cosh( \lambda y)-1}{ y^2}  \right) \langle S\rangle_{k} (y) \right\} \textbf{1}_{\{T =k\}} \right)\nonumber\\
 &\leq& \sum_{k=1}^{n}\mathbf{E}_{\lambda} \left(\exp \left\{-\lambda S_k   + \left(\frac{\cosh( \lambda y)-1}{ y^2} \right) G_k^{y}\right\} \textbf{1}_{\{T =k\}} \right). \label{fgdshu}
\end{eqnarray}
By   the fact $S_{k}\geq x$ and $G_k^{y} \leq v^2$  on the set $\{T(x,y,v) =k\}$, inequality (\ref{fgdshu}) implies that,
for all $x,\lambda, v>0$ and all $y \geq 0$,
\begin{eqnarray}
  && \mathbf{P}( S_k\geq x\ \mbox{and}\  G_k^{y} \leq v^2\ \mbox{for some}\ k\in[1,n]) \nonumber\\
 &\leq&  \exp \left\{-\lambda x   + \left(\frac{\cosh( \lambda y)-1}{ y^2}  \right) v^2 \right\}  \label{lastie1}.
\end{eqnarray}
The last inequality attains its minimum at
\[
\lambda=  \lambda(x)= \frac{1}{y} \log\left( \sqrt{1+\frac{x^2 y^2}{v^4}}+ \frac{x y}{v^2}\right).
\]
Substituting $\lambda = \lambda(x)$ in (\ref{lastie1}), we obtain the desired inequality.
This completes the proof  of Theorem \ref{th3}.\hfill\qed

\section*{Acknowledgements}
We would like to thank Igal Sason  for his helpful remarks and suggestions.

\end{document}